\documentclass[11pt]{article}
\usepackage{enumerate}

\usepackage{epsfig,graphicx,graphics,latexsym,amssymb,amsfonts,amsmath,amscd}
\usepackage{changebar}
\usepackage{graphicx}

\newtheorem{definition}{{\bf Definition}}[section]
\newtheorem{theo}[definition]{{\bf Theorem}}

\newtheorem{coro}[definition]{{\bf Corollary}}
\newtheorem{pro}[definition]{\noindent {\bf Proposition}}

\newtheorem{lem}[definition]{\noindent {\bf Lemma}}

\def\Proof{{\parindent0pt {\bf Proof.\ }}}

\def\endproof{\hfill {\kern 6pt\penalty 500
   \raise -0pt\hbox{\vrule \vbox to5pt {\hrule width 5pt
  \vfill\hrule}\vrule}}}
\begin{document}
\title{ The $(\leq 6)$-half-reconstructibility of digraphs}
\author{Jamel Dammak, Baraa Salem   \\
 { \it  Mathematics Department, Faculty of Sciences of
Sfax}\\
{ \it B.P. 802, 3018 Sfax, Tunisia}\\\\
E-mail : jdammak@yahoo.fr \\
E-mail : salem.baraa@yahoo.fr }



 \maketitle
 \centerline{ \bf Abstract}
 \vskip 0.0000001cm
Let $G=(V,A)$ be a digraph. With every subset $X$ of $V$, we
associate the subdigraph $G[X]=(X,A\cap (X\times X))$
by $X$. Given a positive integer $k$, a digraph $G$ is $(\leq
k)$-half-reconstructible if it is determined up to duality by its
subdigraphs of cardinality $\leq k$.  In 2003,
J. Dammak characterized the $(\leq k)$-half-reconstructible finite digraphs, for $k\in \{7,8,9,10,11\}$.
In this paper, we characterize the $(\leq 6)$-half-reconstructible   digraphs.\\


\section{Introduction}
A {\it directed graph} or simply {\it digraph} $G$ consists
of a set $V(G)$ of vertices together with
a prescribed collection $A(G)$ of ordered pairs of distinct
vertices, called the set of the {\it arcs} of $G$. Such a digraph
is denoted by $(V(G),A(G))$ or simply by $(V,A)$. {\it The
cardinality of} $G$ is that of $V$. We denote this cardinality
by $|V(G)|$ as well as $|G|$. Given a digraph $G=(V,A)$, the \textit{dual}
of $G$ is the digraph $G^*=(V,A^*)$ defined by: for
$x\neq y\in V,\ (x,y)\in A^*$ if $(y,x)\in A$. With each subset $X$ of $V$ is associated the \textit{subdigraph}
$(X,A\cap (X\times X))$ of $G$ \textit{induced} by $X$ denoted by $G[X]$.
The subdigraph $G[V-X]$ is
also denoted by $G-X$. For $x\neq y \in V$,
$x \longrightarrow_{G} y$ or $y \longleftarrow_{G} x$ means
$(x,y)\in A$ and $(y,x)\notin A$, $x \longleftrightarrow_{G} y$
means $(x,y)\in A$ and $(y,x)\in A$, $x\ldots_{G} y$ means
$(x,y)\notin A$ and $(y,x)\notin A$.
 For $x\in V$ and $Y\subseteq V$,
$x \longrightarrow_{G} Y$ signifies that for every $y\in Y$, $x \longrightarrow_{G} y$.
For $X, Y\subseteq V$, $X \longrightarrow_{G} Y$ (or simply $X \longrightarrow Y$ or
$X<Y$ if there is no danger of confusion) signifies that for every $x\in X$, $x \longrightarrow_{G} Y$.
For $x\in V$ and for $X, Y\subseteq V$, $x\longleftarrow_{G} Y$, $x\longleftrightarrow_{G}
Y$, $x\ldots_{G} Y$, $X\longleftrightarrow_{G}Y$ and $X\ldots_{G} Y$
are defined in the same way. Given a digraph $G=(V,A)$, two distinct vertices $x$ and $y$ of $G$ form a \textit{directed pair}
if either $x \longrightarrow_{G} y$ or $x \longleftarrow_{G} y$.
Otherwise, $\{x,y\}$ is a \textit{neutral pair}; it is \textit{full} if $x \longleftrightarrow_{G} y$,
and \textit{void} when $x\ldots_{G} y$.

A digraph $T = (V,A)$ is a {\it tournament} if
all its pairs of vertices are directed.  A \textit{transitive tournament}
is a tournament $T$ such that for $x, y, z \in V(T)$, if $x \longrightarrow_{T} y$
and $y \longrightarrow_{T} z$ then $x \longrightarrow_{T} z$. This is simply a \textit{chain} (that is a set equipped with a linear order  in
which the loops have been deleted). Hence, we will consider the chain $\omega$ of non negative integers as a transitive tournament as well as
the chain $\omega^{*}+\omega$ of integers. A \textit{ flag} is a digraph hemimorphic to $(\{0,1,2\},\{(0,1),(1,2),(2,1)\})$,
a {\it peak} is a digraph hemimorphic to
\- $(\{0,1,2\},\{(0,1),(0,2),(1,2),(2,1)\})$\-  or to $(\{0,1,2\},\{(0,1),(0,2)\})$  (see Figure $1$).
 A {\it diamond} is a digraph hemimorphic to
$(\{0,1,2,3\},\{(0,1),(0,2),(0,3),(1,2),(2,3),(3,1)\})$ (see Figure $2$).

\vspace{- 1.6cm}
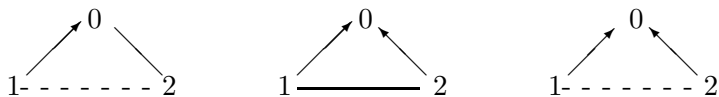
\begin{figure}[htp]\label{peak}
\setlength{\unitlength}{0,9 cm}
\begin{picture}(8,5)
\put(3.5,2){0}
\put(2.3,1){1} \put(4.6,1){2}
\put(2.6,1.3){\vector(1,1){0.8}}
\put(4.6,1.3){\line(-1,1){0.7}}
\put(4.3,1){-}
\put(4,1){-}
\put(3.7,1){-}
\put(3.4,1){-}
\put(3.1,1){-}
\put(2.8,1){-}
\put(2.5,1){-}
\put(7.5,2){0}
\put(6.3,1){1} \put(8.6,1){2}
\put(6.6,1.3){\vector(1,1){0.8}}
\put(8.5,1.3){\vector(-1,1){0.7}}
\put(6.6,1.1){\line(1,0){1.8}}
\put(11.5,2){0}
\put(10.3,1){1} \put(12.6,1){2}
\put(10.6,1.3){\vector(1,1){0.68}}
\put(12.5,1.3){\vector(-1,1){0.7}}
\put(12.3,1){-}
\put(12,1){-}
\put(11.7,1){-}
\put(11.4,1){-}
\put(11.1,1){-}
\put(10.8,1){-}
\put(10.5
,1){-}
\end{picture}
\vspace{- 0.8cm}
\caption{Flag  and peaks.}
\end{figure}

 A {\it prechain}
is a digraph that embeds neither  peak nor diamond nor
 adjacent neutral pairs. Clearly, a chain is a prechain. A \textit{proper prechain} is a prechain that is not a chain.
A prechain which is a tournament is a {\it diamond-free tournament} (that is simply a tournament with no diamond).
Call a \textit{finite consecutivity}
(resp. an \textit{infinite consecutivity}), each digraph on at least three vertices isomorphic to one of the digraphs gotten from a finite chain,
(resp. a transitive tournament of type $\omega$, $\omega^{*}$ or $\omega^{*} +\omega$) such that the pairs of non-consecutive vertices become either all full or all void.
A consecutivity obtained from $\omega$ or  $\omega^{*}$ is called also \textit{one-end infinite consecutivity}.
A {\it cycle} is any digraph isomorphic to one of the digraphs obtained from a finite consecutivity on $n\geq 3$ vertices by
replacing the neutral pair $\{0,n-1\}$ by
$(n-1)\longrightarrow 0$, where $0$ and $n-1$ are the initial and the final extremity respectively. Clearly every $3$-{\it cycle} is
isomorphic to the tournament $(\{0,1,2\},\{(0,1),(1,2),(2,0)\})$.
A \textit{near-chain} is every digraph
obtained from a chain by exchanging the directed pair formed by its
extremities. A $3$-near-chain is a
$3$-consecutivity or a $3$-cycle (see Figure $2$).\\

\begin{figure}[htp]\label{diamond}
\begin{center}
\setlength{\unitlength}{0.8cm}
\begin{picture}(5,5)
\put(-4,2.5){1} \put(1,2.5){2}
\put(-1.6,5.5){0} \put(-1.5,4){3}
\put(1,2.8){\vector(-2,1){2.1}}
\put(-1.4,5.4){\vector(1,-1){2.5}}
\put(-4,2.9){\vector(1,1){2.4}}
\put(-3.8,2.9){\vector(2,1){2.2}}
\put(1,2.7){\vector(-1,0){4.5}}
\put(-1.5,5.4){\vector(0,-1){1}}
\put(4,4){0}
\put(2.8,3){1} \put(5,3){2}
\put(3.1,3.3){\vector(1,1){0.8}}
\put(4.2,4){\vector(1,-1){0.7}}
\put(4.8,3){-}
\put(4.5,3){-}
\put(4.2,3){-}
\put(3.9,3){-}
\put(3.6,3){-}
\put(3.3,3){-}
\put(3
,3){-}
\put(8,4){0}
\put(7,3){2} \put(9,3){1}
\put(7.2,3.2){\vector(1,1){0.8}}
\put(8.2,4){\vector(1,-1){0.8}}
\put(8.8,3){\vector(-1,0){1.5}}
\end{picture}
\vspace{- 2.2cm}
\caption{Diamond, $3$-consecutivity, $3$-cycle.}
\end{center}
\end{figure}
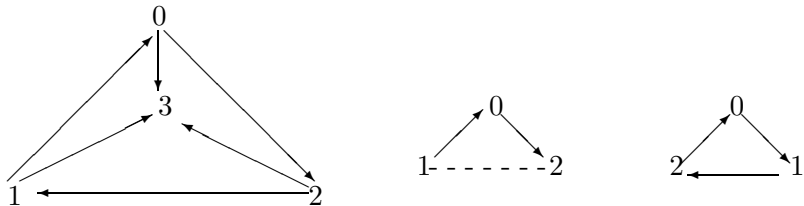

Given a digraph $G=(V,A)$, we define an equivalence relation $\equiv$ on $V$
as follows: for all $x \in V$, $x\equiv x$ and for $x\neq y\in V$,
$x\equiv y$ if there is a sequence $x_{0} = x, . . . , x_{n} = y$ of vertices of $G$
fulfilling that: either $x_i\longrightarrow_{G}x_{i+1}$ or $x_i\longleftarrow_{G}x_{i+1}$,
for all $i\in \{0,...,n-1\}$. If $V\neq \emptyset$, the $\equiv$'s classes are called \textit{arc-connected components} of $G$.
A digraph is said to be \textit{arc-connected} if it has at most one arc-connected component.

Given a digraph $G=(V,A)$, a subset $I$ of $V$
is an \textit{interval} of $G$  if for every $x \in V -I$ either $x\longrightarrow_{G} I$ or $x\longleftarrow_{G} I$ or
$x\longleftrightarrow_{G} I$ or $x\cdots_{G} I$. For instance,
$\emptyset$, $V$ and $\{x\}$ (where $x \in V$)
are intervals of $G$, called {\it trivial intervals}. A digraph is
{\it indecomposable} if all its intervals are trivial, otherwise it is
{\it decomposable}.\\

Given two digraphs $G=(V,A)$ and $G'=(V',A')$, a bijection $f$ from
$V$ onto $V'$ is an {\it isomorphism} from $G$ onto $G'$ provided that for
any $x,y\in V$, $(x,y)\in A$ if and only if
$(f(x),f(y))\in A'$. The digraphs $G$ and $G'$ are then \textit{isomorphic},
which is denoted by $G \simeq G'$, if there
exists an isomorphism from $G$ onto $G'$. If $G$ and $G'$ are not isomorphic, we write $G \not\simeq G'$.
For instance, $G$ and $G'$ are \textit{hemimorphic}, if $G'$ is isomorphic to $G$ or  to $G^*$.
A digraph $G$ is said to be \textit{self-dual} if $G$ is isomorphic to $G^*$.
A digraph $H$ \textit{embeds into} a digraph $G$ if
$H$ is isomorphic to a subdigraph of $G$.\\

Consider two digraphs $G'$ and $G$ on the same vertex set $V$ with $v$ elements and a positive integer $k$.
The digraphs $G'$ and $G$ are $\{k\}$-\textit{hypomorphic} (resp. $\{k\}$-\textit{hemimorphic}) whenever for every subset $X$ of $V$ with $|X| = k$, the subdigraphs $G'[X]$ and $G[X]$ are isomorphic (resp. hemimorphic). $G'$ and $G$ are $\{-k\}$-hypomorphic whenever either $k > v$ or $k\leq v$ and $G'$ and $G$ are $\{v - k\}$-hypomorphic. Notice that $G$ and $G'$ are trivially $\{0\}$-hypomorphic, however $G$ and $G'$ are $\{-0\}$-hypomorphic if and only if they are isomorphic.\\

A digraph $G$ is $\{-k\}$-\textit{self-dual} if
it is $\{-k\}$-\textit{hypomorphic} to $G^*$. Let $F$ be a set of integers. The digraphs $G$ and $G'$ are $F$-\textit{hypomorphic}
(resp. $F$-\textit{hemimorphic}), if for every $k\in F$, the digraphs $G$ and $G'$
are $\{k\}$-hypomorphic (resp.  $\{k\}$-hemimorphic). The digraph $G$ is $F$-reconstructible
(resp. $F$-half-reconstructible) provided that every digraph
$F$-hypomorphic (resp. $F$-hemimorphic) to $G$ is isomorphic (resp. hemimorphic) to $G$.
The digraphs $G$ and $G'$ are $(\leq k)$-\textit{hypomorphic} (resp. $(\leq k)$-hemimorphic) if they are $\{1,...,k\}$-hypomorphic (resp. $\{1,...,k\}$-hemimorphic).
The digraphs $G$ and $G'$ are \textit{hereditarily isomorphic} (resp. \textit{hereditarily hemimorphic})
if for all $X\subseteq V$, $G[X]$ and $G'[X]$ are isomorphic (resp. hemimorphic).\\

Let $G=(V,A)$  and $G'=(V,A')$ be two
$(\leq 2)$-hemimorphic digraphs. Denote  $\mathfrak{D}_{G,G'}$ the
binary relation on $V$ such that: for $x\in V$,
$x \mathfrak{D}_{G,G'}x$; and for $x \neq y\in V$, $x \mathfrak{D}_{G,G'}y$ if
there exists a sequence $x_{0} =x, . . . , x_{n} =y$ of
elements of $V$ satisfying $(x_{i}, x_{i+1})\in A$ if and
only if $(x_{i}, x_{i+1})\notin A'$, for all $i$,
$0 \leq i \leq n - 1$. The relation $\mathfrak{D}_{G,G'}$ is an
equivalence relation called {\it the difference relation},
its classes are called {\it difference classes}.
Let denote $D_{G,G'}$ the set of difference classes.\\

The $(\leq k)$-reconstruction was introduced by R. Fra\"{\i}ss\'{e} in $1970$ \cite{f}. In $1972$, G. Lopez \cite{L1,ls} introduced the difference relation
and showed that:
\begin{theo}\label{thm1} $\cite{L1,ls}$
 The finite digraphs are $(\leq 6)$-reconstructible (i.e: if $G$ and $G'$ are $(\leq 6)$-hypomorphic, then $G$ and $G'$ are isomorphic).
 \end{theo}
In 2002, the $(\leq 5)$-reconstructibility of finite digraphs
was studied by Y. Boudabous \cite{YB}. For $k\in \{3,4\}$, the $(\leq k)$-reconstructibility of finite digraphs
was studied by Y. Boudabous and G. Lopez \cite{YBG} in $2005$.
In $1993$, J. G. Hagendorf raised the $(\leq k)$-half-reconstruction and solved it
with G. Lopez \cite{H}: they proved that, if two digraphs $G$ and $G'$ are $(\leq 12)$-hemimorphic, then either $G'$ and $G$
or $G'$ and $G^{*}$ are $(\leq 6)$-hypomorphic. From that, they obtained in particular: \textit{The finite digraphs are $(\leq 12)$-half-reconstructible}.
Y. Boudabbous and G. Lopez \cite{YBGL} showed that if two finite tournaments $T$ and $T'$
are $(\leq 7)$-hemimorphic, then $T$ and $T'$ are hemimorphic.
Concerning the finite arc-connected digraphs, in $1998$, J. Dammak \cite{JD} proved that they are
$(\leq 7)$-half-reconstructible. He also shown that finite digraphs embedding
a non-self-dual subdigraph of cardinality $k$, are
$(\leq k+6)$-half-reconstructible, for $k\in\{3,4,5,6\}$ \cite{JD}.
M. Pouzet \cite{BD1,BD2} introduced the $\{-k\}$-reconstructibility. P. Ille \cite{P} (resp.
G. Lopez and C. Rauzy \cite{L22}) proved that the finite digraphs on at least $11$ (resp.
$10$) vertices are $\{-5\}$-reconstructible (resp. $\{-4\}$-reconstructible).
Y. Boudabbous \cite{B} improved that: for $k\in \{4,5\}$, two $\{-k\}$-hypomorphic finite tournaments, on at
least $k+6$ vertices, are hereditarily isomorphic.
In $1998$, Y. Boudabbous and J. Dammak \cite{YBJD} introduced the $\{-k\}$-half-reconstruction and proved that:
for $k\in \{4,5\}$, the finite tournaments with at least $k+12$ vertices are $\{-k\}$-half-reconstructible.
In $2012$, Y. Boudabbous and C. Delhomm\'{e} \cite{YBD} studied self duality and introduced the notion of prechain.
In $2003$, J. Dammak \cite{DJC} characterized finite digraphs which are
$(\leq k)$-half-reconstructible, for $k\in \{7,8,9,10,11\}$. After that, N. El Amri
\cite{N}, extended J. Dammak's characterization to infinite digraphs. In the case of tournaments Y. Boudabbous,
A. Boussairi, A. Cha\"{\i}cha\^{a} and N. El Amri \cite{YNB} characterized finite tournaments which are
$(\leq k)$-half-reconstructible, for $k\in \{3,4,5,6\}$.\\

Let $G=(V,A)$ be a digraph and $I$ be a proper interval of $G$. We call  {\it contracted digraph} of $G$ into $I$,
the digraph $G_{I}=((V-I)\cup\{I\},A_{I})$, where $A_{I}$ is defined as follows: $(x,y)\in A_{I}$  if $[( x,y) \in A\cap (V-I)\times (V-I)]$ or $[x=I, y \notin I$ and $\exists z \in I:(z,y)\in A]$
or $[x\notin I, y=I$ and $\exists z \in I:(x,z)\in A]$.
 More precisely, $G_{I}$ is the digraph obtained from $G$ by considering $I$ as a vertex.\\

If $G$ satisfies one of the following conditions, we say that $G$ satisfies the condition $\mathcal{C}_{\infty}$
\begin{enumerate}
    \item[$H_1.$]$G$ has at least an infinite chain interval.
        \item[$H_2.$] $G$ has at least two one-end infinite consecutivity intervals.
        \item [$H_3.$]$G$ has exactly a unique one-end infinite consecutivity interval $I$ and there is no  isomorphism
        $f$ from $G_{I}$ onto $G^*_{I}$ such that $f(I)=I$.
\end{enumerate}
N. El Amri \cite{N} proved that a digraph is non-$(\leq 12)$-half-reconstructible if and only if it verifies $\mathcal{C}_{\infty}$.\\

Given a digraph $G$  with a non-self-dual finite subdigraph,  $\mathcal{C}_{dual}(G)$ denotes the smallest
cardinal of the non-self-dual finite subdigraphs of $G$. From Theorem \ref{thm1},  $3\leq \mathcal{C}_{dual}(G)\leq 6$. In the case where
$G$  has no non-self-dual finite subdigraph,  we set $\mathcal{C}_{dual}(G)=\infty$. \\

Clearly, all non-$(\leq 12)$-half-reconstructible digraphs are not $(\leq k)$-half-reconstructible, for $k\in \{6,7,8,9,10,11\}$.
\begin{theo} $\cite{N}$ Let $G$ be a $(\leq 12)$-half-reconstructible digraph.
The digraph $G$ is non-$(\leq 7)$-half-reconstructible if and only
if one of the following conditions holds:
\begin{enumerate}
\item[$K_1.$] $\mathcal{C}_{dual}(G)=3$ and $G$ admits at least two non-self-dual
arc-connected components which are intervals of type  diamond-free tournament or
non-tournament prechain disjoint from any flag.
\item[$K_2.$] $\mathcal{C}_{dual}(G)=3$ and $G$ has exactly one non-self-dual arc-connected
component $D_0$ which is a diamond-free tournament or
a non-tournament prechain disjoint from any flag, and there is no  isomorphism $f$
from $G_{D_{0}}$ onto $G^*_{D_{0}}$ such that $f(D_{0})=D_{0}$.
\item[$K_3.$] $\mathcal{C}_{dual}(G)=4$ and $G$ admits at least two non-self-dual arc-connected components which are intervals.
\item[$K_4.$] $\mathcal{C}_{dual}(G)=5$ and $G$ admits at least two non-self-dual arc-connected components which are prechain intervals.
\item[$K_5.$] $\mathcal{C}_{dual}(G)=6$ and $G$ admits at least two non-self-dual arc-connected components which are diamond-free tournament  intervals.
\end{enumerate}
 \end{theo}
 As each non-$(\leq 7)$-half-reconstructible digraph is not $(\leq 6)$-half-reconstructible, we obtain our main result:

\begin{theo}\label{T1}  Let $G$ be a $(\leq 7)$-half-reconstructible digraph. The digraph $G$ is non-$(\leq 6)$-half-reconstructible if and only if one of the following conditions holds:
\begin{enumerate}
    \item[$L_1.$] $G$ has at least two intervals $I_{1}$ and $I_{2}$ which
    are non-self-dual  diamond-free tournaments and are not arc-connected components.
    \item[$L_2.$] $G$ has exactly one  non-self-dual interval $I_{0}$ which is a diamond-free tournament that is not an arc-connected component. Furthermore, there is no
 isomorphism $f$ from $G_{I_{0}}$ onto $G^*_{I_{0}}$ such that $f(I_{0})=I_{0}$.
    \item[$L_3.$] $\mathcal{C}_{dual}(G)=3$ and $G$ has at least two non-self-dual arc-connected components $D_1$, $D_{2}$ which are intervals and either disjoint from any
    flag such that $\mathcal{C}_{dual}(G[D_i])=4$,
    for $i\in \{1,2\}$, or these intervals  are non-tournament prechain and each $D_i$ contains a vertex $a_{i}$ of a flag $\{a_{i},b_{i},c_{i}\}$ with $b_{i},c_{i}\notin D_{i}$
and  $\{b_{i},c_{i}\}$ is a directed pair.
    \item[$L_4.$]  $\mathcal{C}_{dual}(G)=3$ and $G$ has exactly one non-self-dual arc-connected component $D_0$ which is an interval being
    either disjoint from any flag such that $\mathcal{C}_{dual}(G[D_0])=4$, or  is
    a non-tournament prechain containing a vertex $a_{0}$ of a flag $\{a_{0},b_{0},c_{0}\}$ with $b_{0}, c_{0}\notin D_{0}$ and
    $\{b_{0},c_{0}\}$ is a directed pair. Furthermore, there is no isomorphism $f$
    from $G_{D_{0}}$ onto $G^*_{D_{0}}$ such that $f(D_{0})=D_{0}$.
\end{enumerate}
\end{theo}
If the condition $L_1$ of Theorem \ref{T1} is satisfied, necessarily $I_{1}$ and $I_{2}$ are disjoint and nontrivial.

Also, in condition $L_3$,  the directed pair $\{b_i,c_i\}$ is disjoint from $D_j$ where $\{j\}=\{1,2\}-\{i\}$.\\

The proof of Theorem  \ref{T1} starts with the case of arc-connected digraph, it is developed in the next section. The general case is treated in section $3$.
\section{The $(\leq 6)$-half-reconstructibility of arc-connected digraphs}

\begin{pro}\label{proposition1}  Let $G$ be a $(\leq 7)$-half-reconstructible arc-connected digraph.
The digraph $G$ is non-$(\leq 6)$-half-reconstructible if and only if $G$
verifies conditions $L_1$ or $L_2$ of  Theorem \ref{T1}.
\end{pro}
Recall some results which will be frequently used in this work.
 \begin{lem}$\cite{L1, L2}$ \label{32} Let $G$ and $G'$ be
two $(\leq 3)$-hypomorphic digraphs, and $C\in D_{G,G'}$.
Then,
\begin{enumerate}
    \item  $G[C]$ is arc-connected and $C$ is an
    interval of $G$ and $G'$.
    \item If $G'[C]\simeq G[C]$ for each $C\in D_{G,G'}$, then $G\simeq G'$.
\end{enumerate}
\end{lem}

\begin{lem}$\cite{YB,YBD,L2}$\label{456} Given an integer
$k\geq4$ and two $(\leq k)$-hypomorphic digraphs $G$ and $G'$, and $C\in D_{G,G'}$, the following assertions hold:
 \begin{enumerate}
    \item If $k=4$, then $G[C]$ is either a consecutivity or cycle
     or a chain or a near-chain or a proper prechain.
    \item If $k=5$, then $G[C]$ is either a consecutivity or cycle
    or a chain or a near-chain or a diamond-free tournament
    or $|C|\leq 6$ and $G[C]$ is a self-dual non-tournament prechain.
     \item If $k=6$, then $G[C]$ is either a consecutivity or cycle
    or a chain or a near-chain or $|C|\leq 7$ and $G[C]$ is
    a self-dual prechain.
    \item If $G[C]$ admits no infinite chain interval,
    then $G'[C]\simeq G^*[C]$.
\end{enumerate}
\end{lem}
From Lemma \ref{456}, we have immediately the following Corollaries.
\begin{coro} \label{454}
Given an integer $k\geq4$, two $(\leq k)$-hypomorphic digraphs
$G$ and $G'$ such that $G$ does not verify $C_\infty$,
and  $C\in D_{G,G'}$  non-self-dual. The following assertions hold:
\begin{enumerate}
    \item If $k=4$, then $G[C]$ is either a one-end infinite consecutivity or a proper
    prechain.
    \item If $k=5$, then $G[C]$ is either a one-end infinite consecutivity
    or a diamond-free tournament.
    \item If $k=6$, then $G[C]$ is a one-end infinite consecutivity.
      \item If $G$ is a prechain,
    then $G$ and $G'$ are hemimorphic.
\end{enumerate}
\end{coro}

Since the equivalence classes of $\mathfrak{D}_{G,G^*}$ are the
arc-connected components of $G$, we have.
\begin{coro} \label{C2} Let $k\geq3$ be an integer and
$G$ be a digraph not verifying the condition $C_\infty$, and
$D$ be a non-self-dual arc-connected component interval of $G$.
\begin{enumerate}
  \item $G[D]$ and $G^*[D]$ are $(\leq k)$-hypomorphic if and
  only if $\mathcal{C}_{dual}(G[D])\geq k+1$.
  \item If $\mathcal{C}_{dual}(G[D]) = 5$, then
  $G[D]$ is a proper
    prechain.
  \item  If $\mathcal{C}_{dual}(G[D])= 6$, then $G[D]$ is a diamond-free tournament.
     \item  If $\mathcal{C}_{dual}(G[D])\geq 7$, then $G[D]$ is
  a one-end infinite consecutivity.
\end{enumerate}
\end{coro}
\begin{pro} $\cite{JD, N}$\label{P1}
Given an integer $k\geq 5$, $G=(V,A)$ and $G'=(V,A')$
 two $(\leq k)$-hemimorphic digraphs and $C\in D_{G,G'}$.
\begin{enumerate}
    \item If $C$ is different from its arc-connected component, then
    $C$ is an interval of $G$ and $G'$ and,
   $G[C]$ and $G'[C]$ are $(\leq k-1)$-hypomorphic.
    \item Let $I_{0}$ be a subset  of $C$ such that
    $|I_{0}|=\mathcal{C}_{dual}(G)$. If $G[I_{0}]$ is non-self-dual,
     then $G'[I_{0}]\simeq G^{\star}[I_{0}]$.
    \item  Given a subset $I_{0}$ of $V$ such that
    $|I_{0}|=\mathcal{C}_{dual}(G)$ and $G[I_{0}]$ is non-self-dual
    such that  $G[I_{0}]\simeq G'[I_{0}]$, we have:
     \begin{enumerate}
      \item  If $G[I_{0}]$ is a flag such that $C\cap I_{0}\neq
        \emptyset$, then $C$ is an interval of $G$ and $G'$ and,
        $G[C]$ and $G'[C]$ are $(\leq k-2)$-hypomorphic.
        \item If $k\geq 6$, then $C$ is an interval of $G$ and $G'$
        and, $G[C]$ and $G'[C]$ are $(\leq h)$-hypomorphic where $h= max(\mathcal{C}_{dual}(G)-1,k-\mathcal{C}_{dual}(G))$.
    \end{enumerate}
\end{enumerate}
\end{pro}
From Proposition \ref{P1}, we obtain the next Corollaries.
\begin{coro} \label{CP} Let $k\geq3$ be an integer,
$G$ and $G'$ be two $(\leq 6)$-hemimorphic digraphs, and
$C\in D_{G,G'}$. $G[C]$ and $G'[C]$ are
$(\leq k)$-hypomorphic if and only if
$\mathcal{C}_{dual}(G[C])\geq k+1$.
\end{coro}
\begin{coro} \label{CP1}
 Let $G$ and $G'$ be two $(\leq 6)$-hemimorphic digraphs
and $C\in D_{G,G'}$. Let $I_0$ be a subset of $V$, such
that $G[I_0]$ is a peak or a flag and
$G'[I_0]\simeq G[I_0]$.
\begin{enumerate}
\item $C$ is an interval of $G$ and $G'$ and $G[C]$ unembed a flag.
\item If $J_0$ is a flag, then $G'[J_0]\simeq G[J_0]$.
\item If $C$ is adjacent at a flag $J_0$, then $G[C]$ and $G'[C]$ are $(\leq 4)$-hypomorphic.
\end{enumerate}
\end{coro}
\Proof
\begin{enumerate}
\item If $C$ is different from its arc-connected component,
Proposition \ref{P1} proves that $C$ is an interval of $G$ and
$G'$ and, $G[C]$ and $G'[C]$ are $(\leq 5)$-hypomorphic.
If $C$ is an arc-connected component, as $G'[I_0]\simeq G[I_0]$, then
from Proposition \ref{P1}, $C$ is an interval of $G$ and $G'$.
In this case, again from Proposition \ref{P1}, if $G[I_{0}]$ is
a flag such that $C\cap I_{0}\neq\emptyset$,
then, $G[C]$ and $G'[C]$ are $(\leq 4)$-hypomorphic, otherwise
$G[C]$ and $G'[C]$ are $(\leq 3)$-hypomorphic.
In consequent, $C$ is an interval of $G$ and
$G'$ and, $G[C]$ and $G'[C]$ are $(\leq 3)$-hypomorphic. Thus,
Corollary \ref{CP} implies that $\mathcal{C}_{dual}(G[C])\geq 4$. So, $G[C]$ unembed a flag.
\item Let $J_{0}=\{a,b,c\}$ such that
 $a\longrightarrow_{G}b$,  $b\longleftrightarrow_{G}c$ and $b\cdots_{G}c$.
 By contradiction, we assume that $G'[J_{0}]\simeq G^*[J_{0}]$.
 So, $b\longrightarrow_{G'}a$. Thus, there exists $C_{0}\in D_{G,G'}$
 such that $a,b\in C_{0}$. As $\{a,b\}$ is not an interval of $G[J_{0}]$, then $J_{0}\subseteq C_{0}$
 which contradicts the first assertion of this corollary.
\item As $C\cap J_{0}\neq\emptyset$ and $G'[J_{0}]\simeq G[J_{0}]$,
from 3.(a) of Proposition \ref{P1}, $G[C]$ and $G'[C]$ are $(\leq 4)$-hypomorphic
\end{enumerate}
\endproof
\begin{coro}\label{cy} Let $G=(V,A)$ and $G'=(V,A')$ be two $(\leq 6)$-hemimorphic
digraphs and $D$ be an arc-connected component of $G$.
Let $I_0\subset V$, such that $\mid I_0\mid=\mathcal{C}_{dual}(G)$, $G[I_0]$ is not self-dual and $G'[I_0]\simeq G[I_0]$.
\begin{enumerate}
\item If $C\in {D}_{G,G'}$, then $C$ is an interval of $G$ and $G'$.
\item If $D$ is not an interval of $G$, then $\mathfrak{D}_{G[D],G'[D]}$  has at least two equivalence classes.
\item If $\mathcal{C}_{dual}(G[D])=3$, then $\mathfrak{D}_{G[D],G'[D]}$ has at least two equivalence classes.
\end{enumerate}
\end{coro}

\Proof
\begin{enumerate}
\item If $C$ is different from its arc-connected component, Proposition \ref{P1} proves that $C$ is an interval of $G$ and $G'$.
If $C$ is an arc-connected component, as $G'[I_0]\simeq G[I_0]$, then
    from Proposition \ref{P1} $C$ is an interval of $G$ and $G'$.
\item If $\mathfrak{D}_{G[D],G'[D]}$ has one equivalence class, $D\in {D}_{G[D],G'[D]}$. As $G'[I_0]\simeq G[I_0]$, then
    from Proposition \ref{P1} $D$ is an interval of $G$ and $G'$.
 \item If $\mathfrak{D}_{G[D],G'[D]}$ has one equivalence class, $D\in {D}_{G[D],G'[D]}$. Since $G'[I_0]\simeq G[I_0]$ and
 $\mid I_0\mid=\mathcal{C}_{dual}(G)=3$, it follows from
 Proposition \ref{P1} that $G[D]$ and $G'[D]$ are $(\leq 3)$-hypomorphic.
  So, from Corollary \ref{CP}, $\mathcal{C}_{dual}(G[D])\geq 4$.
\end{enumerate}
\endproof

\begin{lem} $\cite{JD5, N}$\label{chain}
Let $G=(V,A)$ and $G'=(V,A')$ be two $(\leq 5)$-hemimorphic digraphs.
If $\mathfrak{D}_{G,G'}$ and $\mathfrak{D}_{G^*,G'}$ have just one equivalence class,
then $G$ is a chain.
\end{lem}
\begin{lem}\label{lem2} If a digraph $G$ satisfies $L_1$ or $L_2$, then
$G$ is not $(\leq 6)$-half-reconstructible.
\end{lem}
\Proof
In the two cases, we will construct from $G$ a digraph
$G'$ $(\leq 6)$-hemimorphic and not hemimorphic to $G$.\\
\textbf{Case} $L_1.$  $G$ has at least two non-self-dual intervals $I_1$ and $I_2$ which
 are diamond-free tournaments and not arc-connected components.
    \begin{itemize}
    \item If $G[I_{2}]\simeq G^{\star}[I_{1}]$.
    Let $A=\{I\subset V(G):$ $I$ is an interval of $G$ and $G[I]\simeq G[I_{2}]\}$. For $I_i, I_j \in A$ we have $I_i\cap I_j =\emptyset$.
    Suppose the contrary, as each of the tournaments $G[I_{i}]$ and $G[I_{j}]$ contains at least a $3$-cycle, then $I_i \not\subset I_j$ and $I_j \not\subset I_i$.
    Thus, $I_i- I_j \neq\emptyset$ and  $I_j- I_i \neq\emptyset$ and therefore $I_j- I_i$ and  $I_i\cap I_j$ are two intervals of $G[I_{j}]$, which absurd.
    Let $G'$ be the digraph obtained from $G$ by replacing, every interval
     isomorphic to $G^{\star}[I_{1}]$ by its dual. The digraphs $G$ and
     $G'$ are $(\leq 6)$-hemimorphic. Indeed,
     let $C\in D_{G,G'}$. We have $C$ is an interval of $G$  and $G'$ such that
     $|C|=1$ or, $G[C]\simeq G^{\star}[I_{1}]$ and $G'[C]\simeq G[I_{1}]$.
     As $G[C]$ is a diamond-free tournament, $G[C]$ and $G'[C]$ are
     $(\leq 5)$-hypomorphic. So, for all $X\subset V$ such that $|X|\leq6$,
     if $X\subset C$, $G'[X]\simeq G^{\star}[X]$ otherwise $G'[X]\simeq G[X]$.
    Clearly, $G$ has intervals of type $G[I_{1}]$ and $G^*[I_{1}]$ but
    $G'$ has only intervals of type $G[I_{1}]$, then $G'$ and $G$ are not hemimorphic.
    \item If $G[I_{2}]$ is isomorphic to $G[I_{1}]$, from the first case,
    we may assume that $G$ has no interval isomorphic to $G^{\star}[I_{1}]$. Clearly, the digraph $G'$
    obtained from $G$ by replacing $G[I_{1}]$ by its dual is not hemimorphic to $G$.
    \item If $G[I_{2}]$ and $G[I_{1}]$ are not hemimorphic,
     from the two previous cases, we may suppose that
     $G$ has no interval distinct from $I_1$ and $I_2$ hemimorphic to
     $G[I_{1}]$ or $G[I_{2}]$. The digraph $G'$ obtained from $G$ by replacing $G[I_{1}]$ by
    its dual is not hemimorphic to $G$.
    \end{itemize}
\textbf{Case} $L_2$. $G$ has exactly one non-self-dual interval $I_{0}$
which is a diamond-free tournament and not an arc-connected component
and there is no isomorphism $f$ from $G_{I_{0}}$ onto $G^{\star}_{I_{0}}$
such that $f(I_{0})=I_{0}$.

  Let $G'$ be the digraph obtained from $G$ by replacing $G[I_{0}]$  by
  its dual. As $G[I_{0}]$ is not self-dual, $G'$ is not isomorphic to $G$.
  It suffices to show that $G'$ is not isomorphic to $G^*$, by contradiction let
  $g$ be an isomorphism from $G'$ to $G^*$. Necessarily, $g(I_{0})=I_{0}$.
  So, $g$ induced an isomorphism $f$ from $G_{I_{0}}$ onto
  $G^{\star}_{I_{0}}$ such that $f(I_{0})=I_{0}$ which is absurd.
\endproof

\begin{lem}\label{l3} Let $G=(V,A)$ and $G'=(V,A')$ be two $(\leq 6)$-hemimorphic digraphs such that $G$ does not satisfy the condition $\mathcal{C}_{\infty}$
and $D$ be an arc-connected component of $G$. Let $I_0\subset V$, such that $\mid I_0\mid=\mathcal{C}_{dual}(G)$, $G[I_0]$ is not self-dual and $G'[I_0]\simeq G[I_0]$.
\begin{enumerate}
     \item Let $C\in D_{G[D],G'[D]}$, such that $G[C]$ is neither a one-end infinite consecutivity nor a non-self-dual diamond-free tournament. If $C$ is different from its arc-connected component, then $C$ is an interval of $G$ and $G'$, and $G'[C]\simeq G[C]$.
     \item If G[D] has no interval which is a one-end infinite consecutivity or a non-self-dual diamond-free tournament and if
     $\mathfrak{D}_{G[D],G'[D]}$  has at least two equivalence classes, then for every $C\in D_{G[D],G'[D]}$, $C$ is an interval of $G$ and $G'$,
     and $G'[C]\simeq G[C]$. So, $G'[D]\simeq G[D]$.
\end{enumerate}
\end{lem}
\Proof
\begin{enumerate}
    \item As $C$ is different from its
    arc-connected component, Proposition \ref{P1} proves that $C$
    is an interval of $G$ and $G'$, the subdigraphs $G[C]$ and $G'[C]$ are $(\leq 5)$-hypomorphic.
    As $G[C]$ is neither a one-end infinite consecutivity or a non-self-dual diamond-free tournament, from Corollary \ref{454}, $G'[C]\simeq G[C]$.
 \item Let $C\in D_{G[D],G'[D]}$, as $C$ is different from its arc-connected component, from the first item $C$ is an interval of $G$ and $G'$,
     and $G'[C]\simeq G[C]$. Therefor, from the second assertion of Lemma \ref{32}, $G'[D]\simeq G[D]$.

\end{enumerate}
\endproof

\begin{lem}\label{lem3} Let $G$ and $G'$ be two $(\leq 6)$-hemimorphic
arc-connected digraphs. If $G$ satisfies neither
the condition $\mathcal{C}_{\infty}$ nor $L_1$ nor $L_2$, then $G$ and $G'$ are hemimorphic.
\end{lem}
\Proof
From Lemma \ref{chain},
we may assume that $\mathfrak{D}_{G,G'}$ has at least two classes.
Let $C\in D_{G,G'}$. From Proposition \ref{P1}, $C$ is
an interval of $G$ and $G'$, $G[C]$ and $G'[C]$ are
$(\leq5)$-hypomorphic. As $G$ does not verify $C_{\infty}$,
Lemma \ref{456} proves that $G'[C]\simeq G^*[C]$.
So, if for all $C\in D_{G,G'}$ $C$ is self-dual,
$G'[C]\simeq G[C]$ then
Lemma \ref{32} implies that $G'\simeq G$.
Otherwise,  there exists
a non-self-dual class $C_0\in D_{G,G'}$.
From Corollary \ref{454}, $C_0$ is either
a one-end infinite consecutivity or a diamond-free tournament.
\begin{itemize}
  \item If $C_0$ is a one-end infinite consecutivity,
as $G$ does not verify $C_{\infty}$, $C_0$ is
the unique one-end infinite consecutivity interval of $G$ and there exists
an isomorphism $f$ from $G_{C_0}$ onto $G^{*}_{C_0}$ such
that $f(C_0)=C_0$. Let $C\neq C_0\in D_{G,G'}$.
Clearly $C$ is not a non-self-dual diamond-free tournament; otherwise,
as $G[f(C)]\simeq G^*[C]$, $C$ and $f(C)$ are two non-self-dual intervals of $G$
which are diamond-free tournaments and not arc-connected components
which contradicts the fact that $G$ does not satisfy
the assertion $L_1$.
  \item If $C_0$ is a non-self-dual diamond-free tournament,
  from $L_1$ and $L_2$, $C_0$ is
the unique non-self-dual diamond-free tournament interval of $G$ and there exists
an isomorphism $f$ from $G_{C_0}$ onto $G^{*}_{C_0}$ such
that $f(C_0)=C_0$. Clearly, for all $C\neq C_0\in D_{G,G'}$,
$C$ is not a one-end infinite consecutivity;
otherwise, as $G[f(C)]\simeq G^*[C]$, $C$ and $f(C)$ are two
one-end infinite consecutivity intervals of $G$
which contradicts the fact that $G$ does not satisfy
$C_{\infty}$.
\end{itemize}
From the two cases, for each $C\neq C_0\in D_{G,G'}$,
$C$ is neither a one-end infinite consecutivity
nor a non-self-dual diamond-free tournament.
So, Lemma \ref{l3} proves that
$G'[C]\simeq G[C]$. Therefor, from Lemma \ref{32},
there exists an isomorphism $g$ from $G'_{C_0}$
onto $G_{C_0}$ such that $g(C_0)=C_0$. Thus, $f\circ g$ is an isomorphism
from $G'_{C_0}$ onto $G^*_{C_0}$ such that
$f\circ g(C_0)=C_0$. As $G'[C_0]\simeq G^*[C_0]$, then $G'\simeq G^*$.
\endproof

The proof of Proposition \ref{proposition1} is an immediate consequence of Lemmas \ref{lem2} and \ref{lem3}.
\section{Proof of Theorem \ref{T1}.}

\begin{lem}\label{LP2}   If a digraph $G$ satisfies one of the conditions $L_3$ or $L_4$, then
$G$ is non-$(\leq 6)$-half-reconstructible.
\end{lem}
\Proof
In all these cases, we will construct from $G$ a digraph $G'$ $(\leq 6)$-hemimorphic and not hemimorphic to $G$.\\
  \textbf{Case $L_3.$} $G$
  has at least $G$ has at least two non-self-dual arc-connected components $D_1$, $D_{2}$ which are intervals of $G$ and either disjoint from any
    flag such that $\mathcal{C}_{dual}(G[D_i])=4$,
    for $i\in \{1,2\}$, or non-tournament prechains containing
    a vertex joining two neutral pairs of a flag of $G$.
    \begin{itemize}
    \item If $G[D_{2}]\simeq G^{\star}[D_{1}]$, the digraph  $G'$
    obtained from $G$ by replacing, every  arc-connected component
     isomorphic to $G^{\star}[D_{1}]$ by its dual, is
     $(\leq 6)$-hemimorphic to $G$. Indeed,
     for $C\in D_{G,G'}$,
     $C$ is an interval of $G$ and $G'$ such that
     $|C|=1$ or, $G[C]\simeq G^{\star}[D_{1}]$ and $G'[C]\simeq G[D_{1}]$.
     Thus, for all $X\subset V$ such that $|X|\leq6$,
     if there exists  $C'\in D_{G,G'}$ such that $|X\cap C'|\geq 4$, then $G'[X]\simeq G^{\star}[X]$ otherwise $G'[X]\simeq G[X]$.
     Clearly,  $G$ has intervals of type $G[D_{1}]$ and $G^*[D_{1}]$ but
    $G'$ has only intervals of type $G[D_{1}]$, then $G'$ is not hemimorphic to $G$.

    \item If $G[D_{2}]$ is isomorphic to $G[D_{1}]$, from the first case,
    we may assume that $G$ has no arc-connected component interval isomorphic to $G^{\star}[D_{1}]$. The digraph $G'$
    obtained from $G$ by replacing $G[D_{1}]$ by its dual is not hemimorphic to $G$.
    \item If $G[D_{2}]$ and $G[D_{1}]$ are not hemimorphic,
     from the two previous cases, we may suppose that
     $G$ has no arc-connected component interval distinct from $D_1$ and $D_2$ hemimorphic to
     $G[D_{1}]$ or $G[D_{2}]$. The digraph $G'$ obtained from $G$ by replacing $G[D_{1}]$  by
    its dual is not hemimorphic to $G$.
\end{itemize}
 \textbf{Case $L_4.$}
 $G$ has exactly one non-self-dual arc-connected component $D_0$ which is
    either an interval disjoint from any flag such that $\mathcal{C}_{dual}(G[D_0])=4$ or an interval
    of type non-tournament prechain  containing
    a vertex joining two neutral pairs of a flag of $G$, and there is no  isomorphism $f$
    from $G_{D_{0}}$ to $G^*_{D_{0}}$ such that $f(D_{0})=D_{0}$.

   The digraph $G'$ obtained from $G$ by replacing $G[D_{0}]$  by
   $G^*[D_{0}]$ is not hemimorphic to $G$. Indeed as $G[D_{0}]$ is not self-dual, $G'$ is not isomorphic to $G$.
  It suffices to show that $G'$ is not isomorphic to $G^*$, by contradiction let
  $g$ be an isomorphism from $G'$ to $G^*$. Necessarily, $g(D_{0})=D_{0}$.
  So, $g$ induced an isomorphism $f$ from $G_{D_{0}}$ onto
  $G^{\star}_{D_{0}}$ such that $f(D_{0})=D_{0}$ which is absurd.
\endproof

Conversely, assuming that $G=(V,A)$ does not verify $L_1$, $L_2$, $L_3$, and $L_4$,  we will prove that $G$
 is $(\leq 6)$-half-reconstructible. As $G$ is  $(\leq 7)$-half-reconstructible,
 $G$ is  $(\leq 12)$-half-reconstructible. So, in the sequel, the digraphs considered
 do not satisfy any of the conditions
$C_\infty$, $K_1$, $K_2$, $K_3$, $K_4$, $K_5$, $L_1$, $L_2$, $L_3$, and $L_4$.


\begin{lem}\label{l3200}  If $\mathcal{C}_{dual}(G)\geq 4$,
then $G$ is $(\leq 6)$-half-reconstructible.
\end{lem}
\Proof
Let $G'=(V,A')$ be a digraph $(\leq 6)$-hemimorphic to $G$. Let $D$ be an arc-connected component of $G$. Since $\mathcal{C}_{dual}(G)\geq 4$,
$D$ is an interval of $G$ and $G'$. From Proposition \ref{proposition1},
$G[D]$ and $G'[D]$ are hemimorphic. If $D$ is self-dual, then
$G'[D]\simeq G^*[D]\simeq G[D]$.

$\bullet$ If $\mathcal{C}_{dual}(G)=4$, from $K_3$, $G$  has at most a non-self-dual
arc-connected component.

$\bullet$ If $\mathcal{C}_{dual}(G)=5$, $\mathcal{C}_{dual}(G[D])\geq 5$. So, Corollary \ref{C2} proves that
$G[D]$ is a one-end infinite consecutivity or a proper prechain.
From the assertions $C_{\infty}$ and $K_4$,  $G$ has
at most a non-self-dual arc-connected component.

$\bullet$ If $\mathcal{C}_{dual}(G)=6$, $\mathcal{C}_{dual}(G[D])\geq 6$.
So, Corollary \ref{C2} implies that $D$ is a
a one-end infinite consecutivity or a diamond-free tournament.
As $G$ does not verify none of the assertions $C_{\infty}$ and $K_5$,
$G$ has at most a non-self-dual arc-connected component.

$\bullet$ If $\mathcal{C}_{dual}(G)\geq 7$, then Corollary \ref{C2} proves that $D$ is a
a one-end infinite consecutivity. The condition $C_{\infty}$ proves that
$G$ has at most a non-self-dual arc-connected component.

In consequent, $G$ has at most a non-self-dual arc-connected component $D_0$.
Thus, if $G'[D_0]\simeq G[D_0]$, then $G'\simeq G$
and if $G'[D_0]\simeq G^*[D_0]$, $G'\simeq G^*$.
\endproof \\



\begin{lem}\label{L313} Let $G'=(V,A')$  be a digraph $(\leq 6)$-hemimorphic to $G$
 and $I_0$ be a subset of $V$, such that $G[I_0]$ is a peak or a flag and $G'[I_0]\simeq G[I_0]$. Assume that $G$ has an interval $M_0$
which is either a one-end infinite consecutivity or a non-self-dual diamond-free tournament. Let $D$ be an arc-connected component disjoint from $M_0$. Then,
\begin{enumerate}
\item There exists an isomorphism $f$ from $G_{M_0}$ onto $G^*_{M_0}$ such that $f(M_0)=M_0$.
\item $D$ has not an interval which is either a one-end infinite consecutivity or a non-self-dual diamond-free tournament.
\item $\mathfrak{D}_{G[D],G'[D]}$ has at least two equivalence classes or $D$ is self-dual.
\item $G$ and $G'$ are hemimorphic.
\end{enumerate}
\end{lem}

\Proof
Denote $D_0$ the arc-connected component containing $M_0$. We have $D\neq D_0$.
\begin{enumerate}
  \item
 $\bullet$ If $M_0$ is a one-end infinite
consecutivity interval of $G$,
from $C_{\infty}$, $M_0$ is the unique
one-end infinite consecutivity interval of $G$
and there is an isomorphism $f$ from $G_{M_{0}}$ onto
$G^*_{M_{0}}$ such that $f(M_{0})=M_{0}$.

$\bullet$ If $M_0$ is a non-self-dual diamond-free tournament.

If $M_0$ is not an arc-connected component from $L_1$
and $L_2$, there is an isomorphism $f$ from $G_{M_{0}}$ onto
$G^*_{M_{0}}$ such that $f(M_{0})=M_{0}$.

If $M_0$ is an arc-connected component from $K_1$,
and $K_2$, there is an isomorphism $f$ from $G_{M_{0}}$ onto
$G^*_{M_{0}}$ such that $f(M_{0})=M_{0}$.
 \item By contradiction, assume that $D$ has an interval $I$ which is either a one-end infinite consecutivity or a non-self-dual diamond-free tournament.

 $\bullet$ If $I$ is a one-end infinite consecutivity, as $G[f(I)]\simeq G^*[I]$, $I$ and $f(I)$ are two
one-end infinite consecutivity intervals of $G$ which contradicts the fact that $G$ does not satisfy $C_{\infty}$.

$\bullet$ If $I$ is a non-self-dual diamond-free tournament.

If $I$ is not an arc-connected component, as $G[f(I)]\simeq G^*[I]$, $I$ and $f(I)$ are two non-self-dual intervals of $G$
which are diamond-free tournaments and not arc-connected components
which contradicts the fact that $G$ does not satisfy the assertion $L_1$.

If $I$ is an arc-connected component, as $G[f(I)]\simeq G^*[I]$, $I$ and $f(I)$ are two non-self-dual intervals of $G$
which are diamond-free tournaments and arc-connected components
which contradicts the fact that $G$ does not satisfy the assertion $K_1$.
\item By contradiction, assume that $\mathfrak{D}_{G[D],G'[D]}$ has a unique class and $D$ is non-self-dual.
As $D$ is the unique class of $\mathfrak{D}_{G[D],G'[D]}$, Corollary \ref{cy} implies, $D$ is an interval of $G$ and $\mathcal{C}_{dual}(G[D])\geq 4$.
As $D$ has not an interval which is either a one-end infinite consecutivity or a non-self-dual diamond-free tournament,
then $D$ is not neither a one-end infinite consecutivity or a non-self-dual diamond-free tournament.

$\bullet$ If $D$ is adjacent at a flag $J_0$ of $G$.
Since $\mathcal{C}_{dual}(G[D])\geq 4$ and $D$ is an interval,
$D$ contains only the vertex joining the neutral
pairs of $J_0$. From Corollary \ref{CP1},
$G'[J_0]\simeq G[J_0]$, $G[D]$ and $G'[D]$ are $(\leq 4)$-hypomorphic. As $D\in D_{G,G'}$,
from Corollary \ref{454}, $G[D]$ is a non-tournament prechain. As $G[f(D)]\simeq G^*[D]$,
$D$ and $f(D)$ are two non-self-dual
arc-connected components which are intervals
non-tournament prechain contains a vertex joining two neutral
pairs of a flag of $G$, which contradicts that $G$ does not verify $L_3$.

$\bullet$ If $D$ is disjoint from any flag.

If $\mathcal{C}_{dual}(G[D])=4$, as $G[f(D)]\simeq G^*[D]$,
$D$ and $f(D)$ are two non-self-dual arc-connected
components disjoint from any flag such that
$\mathcal{C}_{dual}(G[D])=C_{n.s}(G[f(D)])=4$, which contradicts
that $G$ does not verify $L_3$.

If $\mathcal{C}_{dual}(G[D])\geq 5$, then, from Corollary \ref{C2},
$G[D]$ is non-tournament prechain disjoint from any flag. As $G[f(D)]\simeq G^*[D]$,
$D$ and $f(D)$ are two non-self-dual arc-connected components which
are intervals non-tournament prechain disjoint from any flag, which contradicts
that $G$ does not verify $K_1$.
\item If $D$ is not an interval, Corollary \ref{cy} proves that
$\mathfrak{D}_{G[D],G'[D]}$ has at least two equivalence classes. If $\mathfrak{D}_{G[D],G'[D]}$  has at least two equivalence classes,
as $D$ has not an interval which is either a one-end infinite consecutivity or a non-self-dual diamond-free tournament,
then, from Lemma \ref{l3}, for each $C\in D_{G[D],G'[D]}$,
$C$ is an interval of $G$ and $G'[C]\simeq G[C]$. If $D$ is an interval self-dual, from Proposition \ref{proposition1},
$G'[D]\simeq G[D]\simeq G^{*}[D]$.
From Lemma \ref{32}, there exists an isomorphism $g$ from
$G'_{D_0}$ onto $G_{D_0}$ such that $g(D_0)=D_0$.
Besides, $f$ induced an isomorphism
  $h$ from $G_{D_0}$ onto $G^*_{D_0}$ such that $h(D_0)=D_0$.
 Further, using Proposition \ref{proposition1}, $G'[D_0]$ and $G[D_0]$ are hemimorphic.
 So, if $G'[D_0]\simeq G[D_0]$,
 the isomorphism $g$ proves that $G'\simeq G$
  and if $G'[D_0]\simeq G^*[D_0]$, the isomorphism $h\circ g$ from
  $G'_{D_0}$ onto $G^*_{D_0}$ such that $h\circ g(D_0)=D_0$ implies that $G'\simeq G^*$.
\end{enumerate}
\endproof
\begin{lem}\label{Lco}
 Let $G'=(V,A')$ be a digraph $(\leq 6)$-hemimorphic to $G$ and $I_0$ be a subset of $V$, such that $G[I_0]$ is a peak or a flag and
 $G'[I_0]\simeq G[I_0]$. If $G$ has a non-self-dual arc-connected component $D_0\in D_{G,G'}$, then
\begin{enumerate}
\item There exists an isomorphism $f$ from $G_{D_0}$ onto
$G^*_{D_0}$ such that $f(D_0)=D_0$.
\item If $D$ is arc-connected component distinct from $D_0$, then $\mathfrak{D}_{G[D],G'[D]}$ has at least two equivalence classes or $D$ is self-dual.
  \item $G$ and $G'$ are hemimorphic.
\end{enumerate}
\end{lem}
\textbf{Proof}
As $D_0\in D_{G,G'}$, Corollary \ref{cy} proves that $D_0$ is interval of $G$ and $G'$.
From Lemma \ref{L313}, we may assume that $G$ has not an interval which is either a one-end infinite consecutivity
or a non-self-dual diamond-free tournament, so $D_0$ is not neither a one-end infinite consecutivity nor a diamond-free tournament.
As $D_0$ is the unique class of $\mathfrak{D}_{G[D_0],G'[D_0]}$, Corollary \ref{cy}  proves that $\mathcal{C}_{dual}(G[D_0])\geq 4$.
\begin{enumerate}
\item $\bullet$ If $D_0$ is adjacent at a flag $J_0$.
  Since $\mathcal{C}_{dual}(G[D_0])\geq 4$ and $D_0$ is an interval of $G$,
 $D_0$ contains only the vertex joining the neutral
 pairs of $J_0$. Then, Corollary \ref{CP1} proves that
$G'[J_0]\simeq G[J_0]$, $G[D_0]$ and $G'[D_0]$ are $(\leq 4)$-hypomorphic. As $D_0\in D_{G,G'}$, Corollary \ref{454}, proves that $D_0$ is
a non-tournament prechain. From $L_3$
and $L_4$, there is an isomorphism
$f$ from $G_{D_{0}}$ onto $G^*_{D_{0}}$ such
that $f(D_{0})=D_{0}$.

$\bullet$ If $D_0$ is disjoint from any flag.

  If $\mathcal{C}_{dual}(G[D_0])=4$, from $L_3$
and $L_4$, there is an isomorphism $f$ from $G_{D_{0}}$ onto
$G^*_{D_{0}}$ such that $f(D_{0})=D_{0}$.

 If $\mathcal{C}_{dual}(G[D_0])\geq 5$, from Corollary \ref{C2},
$D_0$ is a non-tournament prechain. Thus, $K_1$
and $K_2$ implies that there is an isomorphism
$f$ from $G_{D_{0}}$ onto $G^*_{D_{0}}$ such that $f(D_{0})=D_{0}$.
\item The proof is similar to that of item $3$ of lemma \ref{L313}.
\item The proof  is similar to that of item $4$ of lemma \ref{L313}.
\end{enumerate}
\endproof


 \begin{lem}\label{L314} The digraph $G$ is $(\leq 6)$-half-reconstructible.
\end{lem}
\Proof
Let $G'=(V,A')$ be a digraph $(\leq 6)$-hemimorphic to $G$. From Lemma \ref{l3200}, we may
assume that $\mathcal{C}_{dual}(G)=3$. Let $I_0$ be a subset of $V$
such that $G[I_0]$ is a peak or a flag and $G'[I_0]\simeq G[I_0]$.
If $G$ has an interval which is either a one-end infinite consecutivity
or a non-self-dual diamond-free tournament,
from lemma \ref{L313}, $G$ and $G'$ are hemimorphic.
In the sequel, $G$ has no interval which is either a one-end infinite
consecutivity or a non-self-dual diamond-free tournament.
If $G$ has a non-self-dual arc-connected
component $D_0\in D_{G,G'}$, from Lemma \ref{Lco},
$G$ and $G'$ are hemimorphic. Now, we may suppose that
each arc-connected component $D$ of $G$ is self-dual or
$D\not\in D_{G,G'}$.
\begin{itemize}
\item If $D\not\in D_{G,G'}$ or $D$ is not an interval of $G$,
 Corollary \ref{cy} proves that $\mathfrak{D}_{G[D],G'[D]}$
 has at least two equivalence classes, so from Corollary \ref{cy},
for each $C\in D_{G[D],G'[D]}$, $C$ is an interval of $G$ and $G'$. As $G[C]$ is neither a one-end infinite consecutivity or a non-self-dual diamond-free tournament, Lemma \ref{l3} proves $G'[C]\simeq G[C]$.
 \item If $G[D]$ is self-dual and $D$ is an interval of $G$, from Proposition \ref{proposition1},
$G'[D]\simeq G[D]$.
\end{itemize}
 In consequent, from Lemma \ref{32},  $G'$ and $G$ are isomorphic.
\endproof

The proof of Theorem \ref{T1} is an immediate consequence of Lemmas \ref{lem2}, \ref{LP2} and \ref{L314}.

\end{document}